\documentclass{ifacconf}

\usepackage{graphicx}      
\usepackage{natbib}        
\usepackage{amsmath}
\usepackage{amssymb}
\usepackage{units}
\usepackage{color}
\usepackage{prettyref} 
\newrefformat{sec}{\mbox{Section \ref{#1}}}
\newrefformat{tab}{\mbox{Tab. \ref{#1}}}
\newrefformat{fig}{\mbox{Fig. \ref{#1}}}
\newrefformat{rem}{\mbox{Remark \ref{#1}}}
\newrefformat{alg}{\mbox{Algorithm \ref{#1}}}

\usepackage{algorithm}
\usepackage{algpseudocode}
\usepackage{enumerate}
\usepackage{soul,color}
\usepackage{nth}

\newcommand{\T}{\scriptscriptstyle\top}       
\newcommand{\mathmin}{\operatorname*{minimize}}
\newcommand{\mathst}{\text{s.t.}}

\newcommand{\SA}{\text{SA}}
\newcommand{\PAT}{\mbox{PAT}}
\newcommand{\FAT}{\mbox{FAT}}
\newcommand{\PAS}{\mbox{PAS}}
\newcommand{\FAS}{\mbox{FAS}}

\usepackage{dsfont}
\newcommand{\ind}{\mathds{1}}

\usepackage{algorithm}
\usepackage{algpseudocode}

\newtheorem{assumption}{Assumption}
\newtheorem{requirements}{Requirements}
\newtheorem{remark}{Remark}

\usepackage{tikz}
\usetikzlibrary{shapes,arrows}
\usepackage{pgfplots} 
\usepackage{pgfgantt}
\usepackage{pdflscape}
\pgfplotsset{compat=1.5}
\pgfplotsset{plot coordinates/math parser=false}
\newlength\fwidth
\usepackage{grffile}
\usetikzlibrary{plotmarks}
\usetikzlibrary{arrows.meta}
\usepgfplotslibrary{patchplots}

\begin{document}
\begin{frontmatter}

\title{\mbox{\hspace*{-10mm} Nonconvex Consensus ADMM for Cooperative Lane} \mbox{\hspace*{-10mm}Change Maneuvers of Connected Automated Vehicles}}


\author[First]{Alexander Katriniok} 

\vspace*{-1mm}
\address[First]{Ford Research \& Innovation Center (RIC), 
   \mbox{S\"{u}sterfeldstr. 200},\\ 52072 Aachen, Germany (e-mail: de.alexander.katriniok@ieee.org).}%

\begin{abstract}                
Connected and automated vehicles (CAVs) offer huge potential to improve the performance of automated vehicles (AVs) without communication capabilities, especially in situations when the vehicles (or agents) need to be cooperative to accomplish their maneuver. Lane change maneuvers in dense traffic, e.g., are very challenging for non-connected AVs. To alleviate this problem, we propose a holistic distributed lane change control scheme for CAVs which relies on vehicle-to-vehicle communication. The originally centralized optimal control problem is embedded into a consensus-based Alternating Direction Method of Multipliers framework to solve it in a distributed receding horizon fashion. Although agent dynamics render the underlying optimal control problem nonconvex, we propose a problem reformulation that allows to derive convergence guarantees. In the distributed setting, every agent needs to solve a nonlinear program (NLP) locally. To obtain a real-time solution of the local NLPs, we utilize the optimization engine \texttt{OpEn} which implements the proximal averaged Newton method for optimal control (PANOC). Simulation results prove the efficacy and real-time capability of our approach. 
\vspace*{-6mm}
\end{abstract}

\end{frontmatter}

\section{Introduction}
\label{sec:intro}
\vspace*{-2mm}

Automated vehicles (AVs) usually take independent decisions which are based upon sensor measurements and motion predictions of surrounding vehicles.
However, these predictions are often highly uncertain as they rely on simplified assumptions.
This uncertainty may be crucial, especially in situations when the vehicles (hereafter referred to as \textit{agents}) need to rely on these predicted trajectories to accomplish their maneuver. 
For instance, an AV might fail to perform a fully automated lane change or lane merging maneuver when traffic is dense and the target lane is already occupied. When exploiting vehicle-to-vehicle (V2V) communication, we can alleviate this issue by transmitting future control actions or state trajectories to other agents, or even by cooperatively negotiating control actions.

In this paper, we focus on fully automated lane change maneuvers in situations when the target lane is already occupied. During the last two decades, the problem of automating lane change maneuvers has intensively been discussed in literature, see \cite{Bevly2016a} for a comprehensive survey. In the recent years, cooperative control strategies have gained significant attention. Besides those 
based on consensus (\cite{Wang2017a}) or lane change protocols (\cite{An2018a}), optimization based concepts are often a favorable choice as they allow to impose constraints and treat the control problem more holistically. Centralized schemes, which involve a central node that optimizes the agents' control actions, are discussed in \cite{Wang2016} and \cite{Hu2019a}. Decentralized or distributed optimal control schemes, though, may be preferable as they are more resilient and scalable. Existing distributed schemes, such as \cite{Liu2017a,Blasi2018a}, however, introduce conservatism to decouple the agents or require the other agents' state space models to be known.

%
%
%
%
%
%




\vspace*{-1mm}
\subsection{Main Contribution}
\label{sec:intro_contribution}
\vspace*{-3mm}

We propose a distributed 
optimal control approach for collaborative, fully automated lane change maneuvers which adopts the consensus Alternating Direction Method of Multipliers (ADMM) (\cite{Bertsekas1989,Boyd2011}) as methodology to solve the lane change problem in a distributed receding horizon fashion. To 
exchange information, the agents rely on V2V communication. In the considered scenario, the \textit{subject agent} (\SA), i.e., the agent who intends to change lanes determines two consecutive agents in the target lane to eventually merge into the gap between them.
The maneuver is then carried out in two steps: 1) the agents increase their headway distance to allow the SA to change lanes; 2) the SA changes lanes.

Compared to 
the literature, we aim to solve the originally centralized lane change problem in a distributed way without introducing additional conservatism. 
Particularly, every agent optimizes its local control actions while consensus with other agents is achieved through the use of a coordinator, which is run on the SA. As an advantage of our formulation, the agents' parameters and their state space models remain private and do not need to be known by the other agents. This way, we can also reduce load on the communication channel. Moreover, we holistically account for longitudinal and lateral vehicle motion instead of assuming the agents to change instantaneously from one lane to another. Essentially, every agent has to solve a nonconvex nonlinear program (NLP) while the coordinator problem is a standard quadratic program (QP). Although the local NLPs are nonconvex, we propose a problem reformulation that allows to guarantee convergence of the nonconvex consensus ADMM problem. 
We run the algorithm in real-time by adopting \texttt{OpEn} (Optimization Engine) (\cite{Sopasakis2020}) which relies on the proximal averaged Newton method for optimal control (PANOC) (\cite{Stella2017}) to solve the local NLPs fast.

The remainder of the paper is organized as follows. First, we outline the lane change problem in \prettyref{sec:CLCproblem}. Then, we present a centralized formulation in \prettyref{sec:COCP} before we distribute the problem in \prettyref{sec:DOCP}. In \prettyref{sec:results}, we discuss simulation results. Finally, we conclude and give an outlook for future work in \prettyref{sec:conclusion}.

\subsection{Notation}
\label{sec:intro_notation}
With $x_{k+j\mid k}$, we refer to the prediction of variable $x$ at the future time step 
$k+j$ given information up to time $k$ while $x_{\cdot\mid k}$ denotes the trajectory of $x$ 
along the entire prediction horizon of length $N\in\mathbb{N}_{>0}$. For $x\in\mathbb{R}^n$ and $i\in\{1,\ldots, n\}$, $x_i$ is the $i$-th entry of $x$, and the interval $[a,b] \subset \mathbb{N}$ with $a<b$ is denoted as $\mathbb{N}_{[a,b]}$. Finally, $[x]_+ \triangleq \max\{0,x\}$ for $x \in \mathbb{R}$ is referred to as the plus operator.

\section{Cooperative Lane Change Problem}
\label{sec:CLCproblem}

\subsection{Problem Description}
\label{sec:CLCproblem_description}

\begin{figure}[h!]
	\centering
	\def\svgwidth{88mm}
	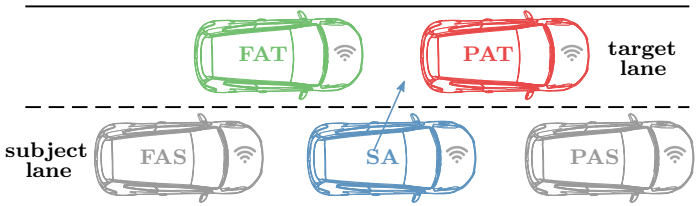	
	\caption{Sketch of the lane change scenario with $N_A=5$ cooperative agents. The SA (blue) aims to merge in between the PAT (red) and FAT (green).}
	\label{fig:CLCproblem_description_sketch}
\end{figure}

\prettyref{fig:CLCproblem_description_sketch} illustrates a sketch of the problem we intend to solve. To reduce complexity, we restrict ourselves to two lane scenarios, in which the \SA{} wants to change lanes while the target lane is already occupied by other agents. This scheme, though, can easily be extended to scenarios with more than two lanes. In our use case, the \SA{} sends a cooperation request to $N_A-1$ agents in order to eventually merge in between two of these agents in the target lane. These are the preceding agent (\PAT{}) and the following agent (\FAT{}) in the target lane in accordance to \prettyref{fig:CLCproblem_description_sketch}. The preceding and following agent in the subject lane are referred to as \PAS{} and \FAS{} respectively. 

Generally, the fully automated lane change maneuver can as such be subdivided into the following steps:
\begin{enumerate}
	\item \textbf{Initialize cooperative group:} SA sets up a cooperative group of $N_A$ agents, in which every agent can communicate with every other agent of the group. 
	\item \textbf{Negotiate agent order:} Determine \PAT{} and \FAT{} \\of \SA{} in the target lane, i.e., the agents that will be in front and behind the SA after changing lanes. 
	\item \textbf{Establish headway:} There must be sufficient headway distance between \PAT{} \& \FAT{} to let the \SA{} in.
	\item \textbf{Conduct lane change:} \SA{} changes lanes.
\end{enumerate}
Some of these steps may also be combined instead of solving them separately. To reduce complexity, we assume that Step 1 and Step 2 have already been accomplished, i.e., the \PAT{} and \FAT{} are known. These steps will further be investigated as part of future work. In this work, we propose a distributed algorithm to solve \mbox{Step 3} and \mbox{Step 4} while relying on the following fundamental assumptions. 
\begin{assumption}\normalfont
	A1.~All agents are equipped with V2V communication; A2.~No communication failures or package dropouts occur; A3.~Agent clocks are synchronized; A4.~Every agent has access to a digital map to have knowledge about the road geometry ahead; A5.~Every agent in the scenario belongs to the cooperative group.
\end{assumption}
\vspace*{-1.0mm}
Assumption A5 reduces complexity in the problem description but does not limit the applicability of our approach. Conversely, it can easily be extended in that direction. 

\subsection{Modeling}
\label{sec:CLCproblem_modeling}

\begin{figure}[b!]
	\centering
	\def\svgwidth{68mm}
	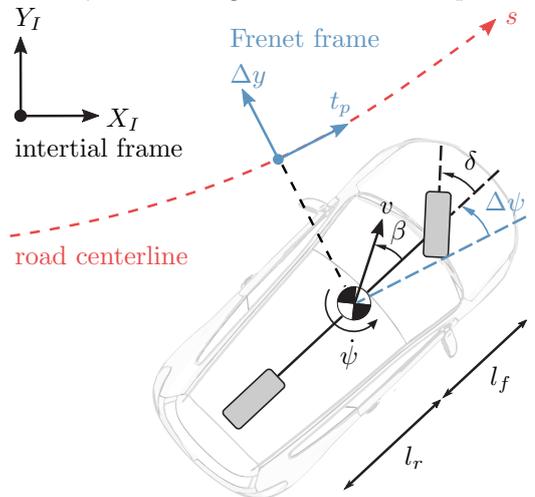	
	\caption{Agent motion in the Frenet frame.}
	\label{fig:modeling_freebody}
\end{figure}
To derive a mathematical model of the lane change maneuver, a kinematic bicycle model (\cite{Rajamani2012}) is adopted to describe the agents' motion in a 
curvilinear reference frame, see \prettyref{fig:modeling_freebody}. For such kind of use case, it is a common approach in literature to apply the Frenet frame (\cite{Qian2016}), in which the agent's position 
is given in terms of its path coordinate $s$ and the perpendicular displacement $\Delta y$ from the road centerline. With this definition, the center of the subject lane in 	\prettyref{fig:CLCproblem_description_sketch} is given by $(s,\Delta y) = (s,-w_\text{lane}/2)$ for any $s$ where $w_\text{lane}$ is the lane width. For every agent, we devise 
a state space model of the form
\begin{subequations}
	\label{eq:modeling_stateSpace}
	\begin{align}
	\hspace*{-1mm}\dot{s} &= v \cos(\Delta \psi+\beta) \left( \frac{1}{1 - \Delta y\, \kappa(s)} \right) \label{eq:modeling_stateSpace_dots}\\
	\hspace*{-1mm}\Delta\dot{y} &= v \sin(\Delta \psi+\beta) \label{eq:modeling_stateSpace_dotdy}\\
	\hspace*{-1mm}\Delta\dot{\psi} &= \frac{v}{l_r} \sin(\beta) - v \cos(\Delta \psi) \left( \frac{\kappa(s)}{1 - \Delta y\, \kappa(s)} \right) \label{eq:modeling_stateSpace_dotdpsi}\\ 
	\hspace*{-1mm}\dot{v} &= a_x \label{eq:modeling_stateSpace_dotv}
	\end{align}
\end{subequations}
where 
\begin{align*}
\beta = \arctan\bigl( \tan(\delta) \, l_r/(l_f+l_r)  \bigr)
\end{align*}
denotes the vehicle sideslip angle, $\delta$ the wheel steering angle, $v$ the vehicle speed, $a_x$ the longitudinal vehicle acceleration, $\Delta \psi$ the heading error between the path tangent (along $t_p$) and the longitudinal vehicle axis. 
Moreover, $l_f$, $l_r$ refer to the distance between the front respectively rear axle and the center of gravity. The path curvature $\kappa(s)$ of the road centerline is assumed to be a known parameterized curve (see assumption A4) in the path \mbox{coordinate $s$.}

The resulting nonlinear state space model 
\(
\dot{x} = f(x,u)
\) features 
the state vector $x=[s,\,\Delta y,\,\Delta\psi,\,v]^{\T}\in \mathbb{R}^{n_x}$ and the input vector  $u=[a_x,\,\delta]^{\T}\in \mathbb{R}^{n_u}$. It can be recognized that model \eqref{eq:modeling_stateSpace} exhibits a singularity for $\Delta y = 1/\kappa(s) = r(s)$ where $r(s)$ is the path radius. In our case, though, even for a radius of $\unit[10]{m}$, a path deviation of $\unit[10]{m}$ is very unlikely to happen. 
To subsequently differ between individual agents, we define the set of collaborative agents
\(
\mathcal{A} \triangleq \{1,\ldots,N_A\}
\)
where $N_A$ is a positive integer. 
This way, we refer to the state vector $x$ of agent $i \in \mathcal{A}$ as $x^{[i]}$.




\section{Centralized Formulation}
\label{sec:COCP}

After generally introducing the cooperative lane change problem in \prettyref{sec:CLCproblem_description}, hereafter, 
we aim to formalize the problem in terms of an optimal control problem (OCP) which is solved in a receding horizon fashion. We start with a centralized formulation 
and derive its distributed variant in \prettyref{sec:DOCP}. Following \prettyref{sec:CLCproblem_description}, we 
subdivide the lane change in two maneuver steps. 
\vspace*{-0.5mm}
\begin{enumerate}[\normalfont(\bfseries M-Step 1\normalfont)]
	\item \SA{}, \PAT{} \& \FAT{} establish required headway.\vspace*{1mm}
	\item The SA changes lanes.
\end{enumerate}
\vspace*{-0.5mm}
We integrate both steps in a single OCP as subsequently outlined.

\subsubsection{Control Objectives}
To perform a proper fully automated cooperative lane change maneuver, the agents have to accommodate certain control objectives. 
First, every agent should track its reference speed $v^{\text{ref}}$ (usually provided by a higher level planning algorithm). Second, the deviation from the lane center, defined through the desired lateral displacement $\Delta y^{\text{ref}}$ from the road centerline, should be minimized. Additionally, we aim to minimize the heading error $\Delta \psi$ (thus choosing $\Delta \psi^{\text{ref}}=0$) to reduce overshooting during a lane change.
For reasons of comfort, we penalize the control input magnitude, that is, the applied longitudinal acceleration $a_x$ and the wheel steering angle $\delta$. 

Along a horizon of $N$ steps, we formalize these objectives for every agent $i \in \mathcal{A}$ as the stage cost at time $k+j$ for $j \in \mathbb{N}_{[0,N-1]}$
\begin{align}
&\ell_{j}^{[i]}({x}_{k+j\mid k}^{[i]},{u}_{k+j\mid k}^{[i]}) 
{}\triangleq{}  {u}_{k+j\mid k}^{[i],\T} \,R^{[i]}\, {u}_{k+j\mid k}^{[i]} \label{eq:problem_localObj_stageCost}  \\[1mm]
&~~~~~~~~~~+  \bigl({x}_{k+j\mid k}^{[i]} - {x}_{k+j\mid k}^{{[i]},\mathrm{ref}} \bigr)^{\T} Q^{[i]}\,  \bigl({x}_{k+j\mid k}^{[i]} - {x}_{k+j\mid k}^{{[i]},\mathrm{ref}} \bigr)
\notag
\end{align}
and the terminal cost
\begin{align}
\ell_{N}^{[i]}({x}_{k+N\mid k}^{[i]}) 
{}\triangleq{} \bigl({x}_{k+N\mid k}^{[i]} - {x}_{k+N\mid k}^{{[i]},\mathrm{ref}} \bigr)^{\T} Q_N^{[i]}\,  \bigl({x}_{k+N\mid k}^{[i]} - {x}_{k+N\mid k}^{{[i]},\mathrm{ref}} \bigr)  \label{eq:problem_localObj_termCost} 
\end{align}
where ${x}_{k+j\mid k}^{[i],\text{ref}}\triangleq[s_{k+j\mid k}^{[i],\text{ref}},\Delta y_{k+j\mid k}^{[i],\text{ref}},\Delta\psi_{k+j\mid k}^{[i],\text{ref}},v_{k+j\mid k}^{[i],\text{ref}}]^{\T}$ 
\mbox{denotes} the reference state while \mbox{$Q^{[i]}\triangleq\mathrm{diag}(q_{s}^{[i]},q_{\Delta y}^{[i]},q_{\Delta \psi}^{[i]},q_{v}^{[i]})$ $\succeq 0$}, 
$Q_N^{[i]}\triangleq\mathrm{diag}(q_{N,s}^{[i]},q_{N,\Delta y}^{[i]},q_{N,\Delta \psi}^{[i]},q_{N,v}^{[i]}) \succeq 0$ and $R^{[i]}\triangleq\mathrm{diag}(r_{a_x}^{[i]},r_{\delta}^{[i]}) \succ 0$ 
are positive (semi)definite weighting matrices. 

\subsubsection{Constraints}
Besides control objectives, we need to ensure that agents only move in their designated lanes. Thus, we constrain 
the lateral displacement $\Delta y$. Moreover, the agents should not exceed the maximum speed $\overline{v}_{k+j\mid k}^{[i]}$ 
at time $k+j$ and not drive backwards. That said, we are able to derive the admissible state set 
\begin{align}
\mathcal{X}_{k+j\mid k}^{[i]} \triangleq  & \left\{ x \in \mathbb{R}^{n_x} ~\lvert~ \Delta \underline{y}_{k+j\mid k}^{[i]} \leq x_2 \leq \Delta \overline{y}_{k+j\mid k}^{[i]} \right. \\ &~~~~~~~~~~~~~~~~~~~\,\land~ \left. 0 \leq x_4 \leq  \overline{v}^{[i]}_{k+j\mid k} \right\} \notag
\end{align}
at time $k+j$ for $j \in \mathbb{N}_{[1,N]}$ where $\Delta\underline{y}_{k+j\mid k}^{[i]}$ and $\Delta\overline{y}_{k+j\mid k}^{[i]}$ are reasonable bounds which assure that agents in neighboring lanes have a minimum lateral distance to each other.
Due to actuator limitations and for reasons of ride comfort and vehicle stability, we constrain the control actions by 
\begin{align}
\mathcal{U}^{[i]} \triangleq \left\{ u \in \mathbb{R}^{n_u} ~\lvert~ \underline{a}_x \leq u_1 \leq \overline{a}_x  ~\land~ \underline{\delta} \leq u_2 \leq \overline{\delta}   \right\}
\end{align} 
where $\underline{a}_x$, $\overline{a}_x$ and $\underline{\delta}$, $\overline{\delta}$ are appropriately chosen 
lower and upper bounds respectively. Moreover, for the same reason, we bound the lateral acceleration (i.e., the product of longitudinal velocity and yaw rate)
\begin{align*}
a_y^{[i]} = (v^{[i]} \cos(\beta^{[i]})) \cdot (v^{[i]} \sin(\beta^{[i]}))/l_r^{[i]}
\end{align*}
and the total acceleration $a_\text{tot}^{[i]}=( (a_x^{[i]})^2 + (a_y^{[i]})^2 )^{\nicefrac{1}{2}}$ through the constraints
\begin{align}
	&-\overline{a}_{y}^{[i]} 
	{}\leq{}
	a_{y,k+j\mid k}^{[i]} 
	{}\leq{}
	\overline{a}_{y}^{[i]}~~\land~~(a_{\text{tot},k+j\mid k}^{[i]})^2 
	{}\leq{} 
	(\overline{a}_{\text{tot}}^{[i]})^2, 
\end{align}
as a function of $(u_{k+j\mid k}^{[i]},x_{k+j\mid k}^{[i]})$, for $j\in \mathbb{N}_{[0,N-1]}$ with appropriate upper bounds $\overline{a}_{y}^{[i]} > 0$ and $\overline{a}_{\text{tot}}^{[i]} > 0$.

Let $\text{Pd}(i) \subset \mathcal{A}$ denote the set of preceding vehicles of Agent $i$. In case Agent $i$ is the \SA{}, $\text{Pd}(i)$ contains the \PAT{} and, until M-Step 2 is completed, the \PAS{}. 
If \mbox{Agent $i$} is the \FAT{}, $\text{Pd}(i)$ contains the \SA{} and, until M-Step 2 is completed, the \PAT{}. For all other agents, $\text{Pd}(i)$ is the agent that is physically ahead. To establish the headway distance between \SA{}, \PAT{} and \FAT{} (M-Step 1) and to avoid collisions between agents, Agent $i$ bounds the headway distance from below, i.e.,
\begin{align}
{s_{k+j\mid k}^{[l]} - s_{k+j\mid k}^{[i]}} \geq d_{\text{hw}},~\forall l \in \text{Pd}(i)
\end{align}
where the lower bound $d_{\text{hw}} > 0$ also encodes the vehicles' geometry. As the agents' path coordinate $s$ refers to the road centerline (see \prettyref{sec:CLCproblem_modeling}), $d_{\text{hw}}$ has additionally to be increased on highly curved road sections in dependence of the (maximum) road curvature $\lvert\kappa(s)\lvert$.


\subsubsection{Optimization Problem}
For notational convenience, hereafter, we stack the states of all agents $i \in \mathcal{A}$ at time $k+j$ in a single vector, i.e.,  
\(
\bar{x}_{k+j\mid k} \triangleq (x_{k+j\mid k}^{[i]})_{i=1}^{N_A}
\). 
Likewise, for the control inputs, we define 
\(
\bar{u}_{k+j\mid k} \triangleq (u_{k+j\mid k}^{[i]})_{i=1}^{N_A}
\) . Additionally, when referring to the entire state or input trajectory, we write
\(
\bar{x}_{\cdot\mid k} \triangleq (\bar{x}_{k+j\mid k})_{j=1}^{N}
\) 
and
\(
\bar{u}_{\cdot\mid k} \triangleq (\bar{u}_{k+j\mid k})_{j=0}^{N-1}
\).
This way, we introduce the aggregated stage costs
\begin{align*}
\ell_{N}(\bar{x}_{k+N\mid k}) &\triangleq \sum_{i=1}^{N_A} \ell_{N}^{[i]}({x}_{k+N\mid k}^{[i]}), \\
\ell_{j}(\bar{x}_{k+j\mid k},\bar{u}_{k+j\mid k}) &\triangleq \sum_{i=1}^{N_A} \ell_{j}^{[i]}({x}_{k+j\mid k}^{[i]},{u}_{k+j\mid k}^{[i]}).
\end{align*}
By adopting a direct multiple shooting formulation (\cite{Bock1984}), that is, including the discrete-time system 
dynamics as equality constraints in the OCP, we phrase the resulting \textbf{centralized lane change NLP} as
\begin{subequations}
	\label{eq:COCP_centralizedNLP}
	\begin{align}
	&\hspace*{-7mm}
	\mathmin_{ \bar{u}_{\cdot \mid k},\,\bar{x}_{\cdot \mid k}  }\ 
	\ell_{N}(\bar{x}_{k+N\mid k}) 
	{}+{} 
	\sum_{j=0}^{N-1} \ell_{j}(\bar{x}_{k+j\mid k},\bar{u}_{k+j\mid k}) 
	\label{eq:COCP_centralizedNLP_cost}
	\\
	\mathst~&~ \text{for every agent~} i \in \mathcal{A}:\notag\\
	&~u_{k+j \mid k}^{[i]} 
	{}\in{} 
	\mathcal{U}^{[i]},~~~~~~~~~~~\,~~~~~~~~~~~~~~~~j\in \mathbb{N}_{[0,N\hspace*{-0.1mm}-\hspace*{-0.1mm}1]}
	\label{eq:COCP_centralizedNLP_inputs}
	\\[0.5mm]
	&~ 
	x_{k+j \mid k}^{[i]} 
	{}\in{} 
	\mathcal{X}_{k+j\mid k}^{[i]},~~~~~~~~~~~\,~~~~~~~~~~~~j\in \mathbb{N}_{[1,N]}
	\label{eq:COCP_centralizedNLP_states}
	\\[0.5mm]
	~&~ x_{k+j+1\mid k}^{[i]} 
	{}={}
	f_d^{[i]}(x_{k+j\mid k}^{[i]},u_{k+j\mid k}^{[i]}),~~~~\,j\in \mathbb{N}_{[0,N-1]} \raisetag{5mm}
	\label{eq:COCP_centralizedNLP_dynamics}  
	\\[0.5mm]		
	&~ 
	-\overline{a}_{y}^{[i]} 
	{}\leq{}
	a_{y,k+j\mid k}^{[i]} 
	{}\leq{}
	\overline{a}_{y}^{[i]},~~~~~~~~~~~~~~\,j\in \mathbb{N}_{[0,N-1]}
	\label{eq:COCP_centralizedNLP_ay}
	\\[0.5mm]
	&~ 
	(a_{\text{tot},k+j\mid k}^{[i]})^2 
	{}\leq{} 
	(\overline{a}_{\text{tot}}^{[i]})^2,~~~~~~~~~~~ ~~~~~j\in \mathbb{N}_{[0,N-1]}
	\label{eq:COCP_centralizedNLP_ahor}
	\\[0.5mm]
	&~ {s_{k+j\mid k}^{[l]}  - s_{k+j\mid k}^{[i]}} \geq d_{\text{hw}},\,\forall l \in \mathrm{Pd}(i); \,j\in \mathbb{N}_{[1,N]} \label{eq:COCP_centralizedNLP_CA} \\[0.5mm]
	&~
	x_{k\mid k}^{[i]} 
	{}={} 
	x_{k}^{[i]},
	\end{align}
\end{subequations}
where 
\begin{align}
	f_d^{[i]}(x_{k+j\mid k}^{[i]},u_{k+j\mid k}^{[i]}) \triangleq x_{k+j\mid k}^{[i]} + \int_{t_{k+j}}^{t_{k+j+1}} \hspace*{-2.5mm}f(x^{[i]}(\tau),u_{k+j\mid k}^{[i]}) \,d\tau \raisetag{3mm} 
	\label{eq:COCP_discreteSysDyn}
\end{align}
represents the discretized system dynamics using zero order hold. We approximate the integral in \eqref{eq:COCP_discreteSysDyn} using a 4th order Runge-Kutta method (\cite{Nocedal2006}).

\begin{remark} \normalfont
During \mbox{M-Step 1}, the SA's reference value $\Delta y^{\text{ref}}$ is set to the center of the subject lane while the bounds on $\Delta y$ prevent the SA from leaving its lane. When the required headway distance between the \SA{}, \PAT{} and \FAT{} has been established (for time $t_k$ and the entire prediction horizon), 
the bounds on $\Delta y$ are modified in \mbox{M-Step 2} to allow the SA to drive in both lanes. Simultaneously, the reference $\Delta y^{\text{ref}}$ is set to the center of the target lane to initiate the lane change. 
After M-Step 2, the bounds on $\Delta y$ are finally adapted to the target lane.
\end{remark}
\vspace*{-0.5mm}
\section{Distributed Solution using ADMM}
\label{sec:DOCP}

Compared to a centralized formulation, as described in \prettyref{sec:COCP}, distributing computations among the agents appears to be more scalable and resilient.
For a distributed control scheme that could eventually be implemented in a test vehicle, we impose the following requirements.
\begin{requirements}\normalfont
	{R1.}~The agents' state space models and their parameters should be private; {R2.}~ Instead of the entire state vector, only position information should be exchanged to reduce communication load; {R3.}~ The distributed OCP must be solved within the sampling time, including the time for data exchange via V2V.
\end{requirements}
\vspace*{-0.5mm}
To solve Problem \eqref{eq:COCP_centralizedNLP} in a distributed fashion, we apply the consensus Alternating Direction Method of Multipliers (ADMM) (\cite{Bertsekas1989}). 
As shown in the remainder of this section, we decompose Problem \eqref{eq:COCP_centralizedNLP} such that objectives and constraints, belonging to the individual agent, 
are incorporated in a local OCP while the joint satisfaction of the minimum headway distance 
is tackled by a {coordinator} which is run on the SA. 


\subsection{Problem Reformulation and Decomposition}
\label{sec:DOCP_decomp}

To come up with a distributed OCP that is compliant with a consensus ADMM formulation, in a first step, we need to reformulate Problem
\eqref{eq:COCP_centralizedNLP}. For notational convenience, we define an augmented optimization vector that contains control input and state trajectories over the prediction horizon, that is, 
\mbox{$\xi_{\cdot\mid k}^{[i]}= (x_{k+j\mid k}^{[i]},u_{k+j-1\mid k}^{[i]})_{j=1}^{N}$}. Moreover, the aggregated vector
\( 
\bar{\xi}_{\cdot\mid k} \triangleq ({\xi}_{\cdot\mid k}^{[i]})_{i=1}^{N_A}
\)
stacks the variables $\xi_{\cdot\mid k}^{[i]}$ 
for every agent $i \in \mathcal{A}$ in a single vector. This way, we can summarize equality constraints \eqref{eq:COCP_centralizedNLP_dynamics} in a compact form as \mbox{$h^{[i]}({\xi}_{\cdot\mid k}^{[i]}) = 0$}. Likewise, inequality constraints \eqref{eq:COCP_centralizedNLP_ay} and \eqref{eq:COCP_centralizedNLP_ahor} can concisely be stated as $g^{[i]}({\xi}_{\cdot\mid k}^{[i]}) \leq 0$. That said, we introduce the associated indicator functions
\begin{align}
\ind_{h^{[i]}}({\xi}_{\cdot\mid k}^{[i]}) &\triangleq \begin{cases}
0, & h^{[i]}({\xi}_{\cdot\mid k}^{[i]}) = 0\\
\infty, & \text{otherwise},
\end{cases} \label{sec:DOCP_decomp_indEqCons}\\
\ind_{g^{[i]}}({\xi}_{\cdot\mid k}^{[i]})
&\triangleq \begin{cases}
0, & g^{[i]}({\xi}_{\cdot\mid k}^{[i]}) \leq 0\\
\infty, & \text{otherwise}. \label{eq:DOCP_decomp_indIneqCons}
\end{cases} 
\end{align}
To accommodate equality and inequality constraints, we define the augment cost for every agent $i$ as
\begin{align}
\phi^{[i]}({\xi}_{\cdot\mid k}^{[i]}) {}\triangleq{}     
\ell_{N}^{[i]}(x_{k+N\mid k}^{[i]}) 
{}&+{} 
\sum_{j=0}^{N-1} \ell_{j}^{[i]}(x_{k+j\mid k}^{[i]},{u}_{k+j\mid k}^{[i]})
\label{eq:DOCP_decomp_augAgentCost}
\\ 
{}&+{}
\ind_{h^{[i]}}({\xi}_{\cdot\mid k}^{[i]}) + \ind_{g^{[i]}}({\xi}_{\cdot\mid k}^{[i]}). \notag
\end{align}
Moreover, we define a cost associated with the minimum headway distance constraints \eqref{eq:COCP_centralizedNLP_CA}, i.e.,
\begin{align}
\gamma(\bar{\xi}_{\cdot\mid k}) \triangleq \sum_{j=1}^{N} \sum_{i=1}^{N_A} \sum_{l \in \mathrm{Pd}(i)} \ind_\text{hw}(\xi_{k+j \mid k}^{[i]},\xi_{k+j \mid k}^{[l]})
\label{eq:DOCP_decomp_hwConsCost}
\end{align}
with the indicator function 
\begin{align*}
\ind_\text{hw}(\xi_{k+j \mid k}^{[i]},\xi_{k+j \mid k}^{[l]})
&\triangleq \begin{cases}
0, & {s_{k+j\mid k}^{[l]}  - s_{k+j\mid k}^{[i]}} \geq d_{\text{hw}} \\
\infty, & \text{otherwise}.
\end{cases} 
\end{align*}
As the minimum headway distance constraints are convex, so is the associated indicator function $\ind_\text{hw}$ and as such the cost $\gamma$. With \eqref{eq:DOCP_decomp_augAgentCost} and \eqref{eq:DOCP_decomp_hwConsCost}, we can rewrite the centralized NLP \eqref{eq:COCP_centralizedNLP} as a \textbf{box constrained NLP}
\begin{subequations}
	\label{eq:DOCP_decomp_boxConsCNLP}
	\begin{align}
	\mathmin_{ \bar{\xi}_{\cdot \mid k} }\ & 
	\gamma(\bar{\xi}_{\cdot\mid k}) + \sum_{i=1}^{N_A}\phi^{[i]}(\xi_{\cdot\mid k}^{[i]}) 
	\label{eq:DOCP_decomp_boxConsCNLP_cost}
	\\
	\mathst~&~ \xi_{k+j \mid k}^{[i]} 
	{}\in{} 
	{\Xi}_{k+j\mid k}^{[i]},~~~~~~~~ i \in \mathcal{A},~ j \in \mathbb{N}_{[1,N]} \label{eq:DOCP_decomp_boxConsCNLP_xiCons}
	\end{align}
\end{subequations}
where the original input and state constraints \eqref{eq:COCP_centralizedNLP_inputs} and \eqref{eq:COCP_centralizedNLP_states} are represented by the close and convex feasible set
\begin{align*}
{\Xi}_{k+j\mid k}^{[i]} \hspace*{-0.3mm}\triangleq\hspace*{-0.3mm} \Bigl\{ (x_{k+j\mid k}^{[i]},u_{k+j-1\mid k}^{[i]}) \hspace*{-0.4mm} \in \hspace*{-0.4mm} (\mathcal{X}_{k+j\mid k}^{[i]} \hspace*{-0.6mm} \times \hspace*{-0.5mm} \mathcal{U}^{[i]} ) \subseteq \mathbb{R}^{n_x} \hspace*{-1.1mm} \times \hspace*{-0.9mm}\mathbb{R}^{n_u} \hspace*{-0.5mm}\Bigr\}.
\end{align*}

To decouple the agents, we introduce new auxiliary variables \mbox{$z_{\cdot\mid k}^{[i]} \triangleq (z_{k+j\mid k}^{[i]})_{j=1}^N$} respectively $\bar{z}_{\cdot\mid k} \triangleq (z_{\cdot\mid k}^{[i]})_{i=1}^{N_A}$ which correspond to Agent $i$'s path coordinate $s^{[i]}$ over the prediction horizon. By imposing the \textit{consensus constraint} $c_{\text{cc}}^{\T} \,\xi_{k+j\mid k}^{[i]} = z_{k+j\mid k}^{[i]}$ with \mbox{$c_{\text{cc}} \triangleq [1~0~0~0~0~0]^{\T}$} for every \mbox{$j \in \mathbb{N}_{[1,N]}$} and every agent $i \in \mathcal{A}$, we can rewrite \eqref{eq:DOCP_decomp_boxConsCNLP} as an \textbf{equivalent problem with auxiliary variables}
\begin{subequations}
	\label{eq:DOCP_decomp_decoupNLP}
	\begin{align}
	\mathmin_{ \bar{\xi}_{\cdot \mid k},\,\bar{z}_{\cdot \mid k}}\ &
	\gamma(\bar{z}_{\cdot\mid k}) + \sum_{i=1}^{N_A}\phi^{[i]}(\xi_{\cdot\mid k}^{[i]})
	\label{eq:DOCP_decomp_decoupNLP_cost}
	\\
	\mathst~&~ c_{\text{cc}}^{\T} \,\xi_{k+j\mid k}^{[i]} = z_{k+j\mid k}^{[i]},~~~~ i \in \mathcal{A},~ j \in \mathbb{N}_{[1,N]} \label{eq:DOCP_decomp_decoupNLP_consensus}\\
	&~ \xi_{k+j \mid k}^{[i]} 
	{}\in{} 
	{\Xi}_{k+j\mid k}^{[i]},~~~~~~~~ i \in \mathcal{A},~ j \in \mathbb{N}_{[1,N]} \label{eq:DOCP_decomp_decoupNLP_xiCons}
	\end{align}
\end{subequations}
By virtue of Problem \eqref{eq:DOCP_decomp_decoupNLP}, it can be recognized that the agents' cost $\phi^{[i]}$ solely depends on the local optimization variable $\xi^{[i]}_{\cdot\mid k}$. Conversely, the cost $\gamma$, which accommodates the minimum headway distance constraints, relies on the auxiliary variables $z_{\cdot\mid k}^{[i]}$ or  \textit{copies} of the agents' path coordinate $s^{[i]}$ over the prediction horizon. Such formulation motivates a distributed solution of Problem \eqref{eq:DOCP_decomp_decoupNLP}, that is, every agents optimizes its local OCP while the minimum headway distance constraints are accommodated by a coordinator. 
For reasons of brevity, hereafter, we abbreviate constraint \eqref{eq:DOCP_decomp_decoupNLP_consensus}  as $C_{\text{cc}} \xi_{\cdot\mid k}^{[i]} = z_{\cdot\mid k}^{[i]}$  where $C_{\text{cc}}$ is a matrix of appropriate dimension.

\vspace{0.5mm}
\begin{remark}\normalfont
Our choice of $c_{\text{cc}}$ in \eqref{eq:DOCP_decomp_decoupNLP_consensus} is crucial to satisfy requirements R1 and R2. This way, only the agents' path coordinates $s^{[i]}=c_{\text{cc}}^{\T}\,\xi^{[i]}$ need to be exchanged amongst each other via V2V. 
\end{remark}

\subsection{Consensus ADMM Framework}
\label{sec:DOCP_ADMM}


To embed the reformulated Problem \eqref{eq:DOCP_decomp_decoupNLP} in the ADMM framework, we \textit{dualize} consensus constraint \eqref{eq:DOCP_decomp_decoupNLP_consensus} and obtain the Augmented Lagrangian function (\cite{Bertsekas1989})
\begin{align} 
&\mathcal{L}_\rho(\bar\xi_{\cdot\mid k},\bar{z}_{\cdot\mid k},\bar{\lambda}) = \gamma(\bar{z}_{\cdot\mid k}) + \sum_{i=1}^{N_A} \phi^{[i]}(\xi_{\cdot\mid k}^{[i]})   \label{eq:DOCP_ADMM_augLagrangian}\\
&~~~~+ \sum_{i=1}^{N_A} \Bigl[ (\lambda^{[i]})^{\T} ({C}_{\text{cc}} \xi_{\cdot\mid k}^{[i]}-z_{\cdot\mid k}^{[i]}) 
+ \frac{\rho}{2} \| {C}_{\text{cc}} \xi_{\cdot\mid k}^{[i]}-z_{\cdot\mid k}^{[i]} \|^2 \Bigr]\notag
\end{align}
where $\rho > 0$ is a constant penalty parameter, $\lambda^{[i]}$ is the vector of Lagrangian multipliers associated with consensus constraint \eqref{eq:DOCP_decomp_decoupNLP_consensus}, and $\bar\lambda \triangleq ( \lambda^{[i]} )_{i=1}^{N_A}$ is the aggregated vector of multipliers.
The consensus ADMM algorithm, applied to minimize \eqref{eq:DOCP_ADMM_augLagrangian} subject to 
\(
\xi_{k+j \mid k}^{[i]} 
{}\in{} 
{\Xi}_{k+j\mid k}^{[i]}
\) 
for \mbox{$i \in \mathcal{A}$} and $j \in \mathbb{N}_{[1,N]}$, 
is summarized in \prettyref{alg:DOCP_ADMM_algo}. In a receding horizon fashion, we run the following steps at time $k$.

\textbf{Step 1:}
Starting with an initial guess 
\( (\xi_{\cdot\mid k}^{[i]}, {z}_{\cdot\mid k}^{[i]}, \lambda^{[i]})  \), 
every agent solves the box constrained NLP \eqref{eq:DOCP_ADMM_algo_agentNLP} and transmits the optimized path coordinate trajectory $s_{\cdot\mid k}^{[i]} = {C}_{cc}\, \xi_{\cdot\mid k}^{[i]}$ to the coordinator. 
As $\phi^{[i]}$ in \prettyref{eq:DOCP_ADMM_augLagrangian} is nonsmooth and nonconvex, NLP \eqref{eq:DOCP_ADMM_algo_agentNLP} adopts a smooth but still nonconvex reformulation $\Phi^{[i]}$ to guarantee convergence of the ADMM scheme, see \prettyref{sec:DOCP_PANOC} and \prettyref{sec:DOCP_convergence}.

\textbf{Step 2:}
The coordinator (residing on the \SA{}) solves the constrained coordination QP \eqref{eq:DOCP_ADMM_algo_coordQP} which imposes the minimum headway distance (encoded in $\gamma$) as linear constraints (cf. \eqref{eq:COCP_centralizedNLP_CA}). For notational convenience, we abbreviate this constraint as $\bar{C}_{\text{hw}} \bar{z}_{\cdot\mid k} \leq \bar{d}_{\text{hw}}$ where $\bar{C}_{\text{hw}}$ and $\bar{d}_{\text{hw}}$ are a matrix and a vector of appropriate dimension. Mostly, the initial condition of the lane change maneuver may look like in \prettyref{fig:CLCproblem_description_sketch}. Then, the original constraint \eqref{eq:COCP_centralizedNLP_CA} may not be satisfied. As a consequence, we may not be able to establish consensus without either violating agent dynamics (as we would need to shift agents' positions instantaneously to satisfy \eqref{eq:COCP_centralizedNLP_CA}) or the minimum headway distance constraints. For this reason, we reformulate the minimum headway distance constraint as a soft constraint 
\begin{align}
\bar{C}_{\text{hw}} \bar{z}_{\cdot\mid k} \leq \bar{d}_{\text{hw}} + \bar{E}_{\text{hw}} \epsilon~~\text{with}~~\epsilon \triangleq (\epsilon_{i,l})_{i\in\mathcal{A},l\in\mathrm{Pd}(i)}
\label{eq:DOCP_ADMM_hwDistConsSoft}
\end{align}
where $\epsilon$ is a vector of slack variables $\epsilon_{i,l} \geq 0$ for every $(i,l)$ with $i \in \mathcal{A}$ and $l \in \mathrm{Pd}(i)$ and $\bar{E}_{\text{hw}}$ is a matrix of appropriate dimension. 
That way, consensus can be established and $\epsilon_{i,l} \rightarrow 0$ will hold after \mbox{M-Step 1} is completed. Moreover, we augment the cost function in QP \eqref{eq:DOCP_ADMM_algo_coordQP} with the additional linear cost term $q_\epsilon^{\T}\epsilon \geq 0$ that penalizes $\epsilon$ where $q_\epsilon$ is a weighting vector of appropriate dimension with all weights larger than zero.
The resulting Problem \eqref{eq:DOCP_ADMM_algo_coordQP} is a standard QP that optimizes over $(\bar{z}_{\cdot\mid k},\epsilon)$. Its solution can be obtained fast using mature QP solvers. 

\begin{algorithm}[t!]
	\caption{Nonconvex ADMM Problem at time $k$}
	\label{alg:DOCP_ADMM_algo}
	\begin{algorithmic}
		\State Initial guess:~~ \( \bar\xi_{\cdot\mid k}, \bar{z}_{\cdot\mid k}, \bar\lambda \) (every agent \& coordinator) 
		\Repeat \vspace*{1mm}
		\State 1) \textbf{Every agent $i \in \mathcal{A}$:} Solve NLP in parallel \vspace*{1mm}
		\begin{align} ~~{\xi}_{\cdot\mid k}^{[i]} \gets \text{arg\,}  \underset{\xi_{\cdot\mid k}^{[i]}}{\text{min}}~ 
		\Phi^{[i]}(\xi_{\cdot\mid k}^{[i]};\alpha,\mu) &+ (\lambda^{[i]})^{\T} {C}_{\text{cc}}\, \xi_{\cdot\mid k}^{[i]} \label{eq:DOCP_ADMM_algo_agentNLP} \\[-3mm] &+ \frac{\rho}{2} \| {C}_{\text{cc}}\, \xi_{\cdot\mid k}^{[i]}-z_{\cdot\mid k}^{[i]} \|^2\notag	\end{align}\\\vspace*{1mm}
		~~~~~~~~~~~~~~~~~~~~~~~\mathst ~~~ \({\xi}_{\cdot \mid k}^{[i]} \in  {\Xi}_{\cdot\mid k}^{[i]} \) \\\vspace*{1mm}
		~~~~~~~~~and transmit ${C}_{\text{cc}}\, \xi_{\cdot\mid k}^{[i]}$ to coordinator.\\
		\State 2) \textbf{Coordinator:} Solve coordination QP
		\begin{align} \bar{z}_{\cdot\mid k} \gets \text{arg\,} \underset{ \bar{z}_{\cdot\mid k}, \epsilon  }{\text{min}}~ 
		\sum_{i=1}^{N_A} \Bigl[  &-(\lambda^{[i]})^{\T} z_{\cdot\mid k}^{[i]} \label{eq:DOCP_ADMM_algo_coordQP}\\&+ \frac{\rho}{2} \| {C}_{\text{cc}}\, \xi_{\cdot\mid k}^{[i]}-z_{\cdot\mid k}^{[i]} \|^2 \Bigr] + q_\epsilon^{\T}\epsilon\notag
		\end{align} \\ \vspace{1mm} 
		~~~~~~~~~~~~~~~~~~~~~~~\mathst ~~~ \( \bar{C}_{\text{hw}} \bar{z}_{\cdot\mid k} \leq \bar{d}_{\text{hw}}  + \bar{E}_{\text{hw}} \epsilon,~\epsilon_{i,l} \geq 0.\) \vspace*{2mm}
		\State 3) \textbf{Coordinator:} Perform dual gradient step
		\begin{align}
		\lambda^{[i]} \gets \lambda^{[i]} + \rho \, ({C}_{\text{cc}}\, \xi_{\cdot\mid k}^{[i]}-z_{\cdot\mid k}^{[i]}),~\forall i \in \mathcal{A}.~\,\label{eq:DOCP_ADMM_algo_coordDualGradStep}\end{align}
		\State 4) \textbf{Coordinator:} Broadcast \( (\bar\lambda, \bar{z}_{\cdot\mid k}) \) to every agent.\vspace{2mm}
		\Until{stopping criterion satisfied}
	\end{algorithmic}
\end{algorithm}
\textbf{Step 3 \& 4:}
The coordinator updates the dual variables $\lambda^{[i]}$ and transmits the vectors \( (\bar\lambda, \bar{z}_{\cdot\mid k}) \) to the agents. 

This scheme is iterated until the stopping criteria, adopted from \cite[Sec. 3.1.1]{Boyd2011}, is satisfied, that is, until the norm of the primal and dual residuals are below their 
thresholds $\epsilon^{\text{prim}} \geq 0$ and $\epsilon^{\text{dual}} \geq 0$ respectively
\begin{align}
\hspace*{-2.3mm}
\| \xi_{\cdot\mid k} - \bar{z}_{\cdot\mid k} \| \leq  \epsilon^{\text{prim}},\,
\| ( \rho \, [{C}_{\text{cc}}\, {\xi}_{\cdot\mid k}^{[i]} - {z}_{\cdot\mid k}^{[i]}] )_{i=1}^{N_A} \| \leq  \epsilon^{\text{dual}}. 
\end{align}
After convergence, every agent $i$ applies the control action $u_{k\mid k}^{[i]\,\star}$ locally. At the next time step $k+1$, Algorithm 1 is warm-started by exploiting the solution $( \bar\xi_{\cdot\mid k}^\star, \bar{z}_{\cdot\mid k}^\star, \bar\lambda^\star )$ from time step $k$ as initial guess. 

\subsection{Agent NLPs: Smooth Reformulation and Solution}
\label{sec:DOCP_PANOC}

To ensure convergence of \prettyref{alg:DOCP_ADMM_algo}, $\phi^{[i]}$ in \prettyref{eq:DOCP_ADMM_augLagrangian} needs to be convex or at least a smooth nonconvex function, see \cite{Hong2016}. Convexity, though, is not satisfied due to nonconvex system dynamics $h^{[i]}(\xi_{\cdot\mid k+j}^{[i]})=0$. At the same time, indicator functions in \eqref{eq:DOCP_decomp_augAgentCost} are nonsmooth. Therefore, we come up with a smooth reformulation of the indicator functions $\ind_{h^{[i]}}$ and $\ind_{g^{[i]}}$.

By applying the Augmented Lagrangian method (ALM) \cite[Chap.~17]{Nocedal2006}, we replace the indicator function $\ind_{h^{[i]}}({\xi}_{\cdot\mid k}^{[i]})$ of the equality constraints, related to system dynamics, with the cost
\begin{align}
	\phi_\text{ALM}^{[i]}({\xi}_{\cdot\mid k}^{[i]};\alpha,\mu) {}\triangleq{} \mu^{\T} h^{[i]}({\xi}_{\cdot\mid k}^{[i]}) +
	\frac{\alpha}{2} \,\| h^{[i]}({\xi}_{\cdot\mid k}^{[i]}) \|^2
	\label{eq:DOCP_ADMM_PANOC_ALMcons}
\end{align}
where $\alpha > 0$ is a penalty parameter and $\mu$ the vector of Lagrangian multipliers related to the equality constraints. 
Moreover, we rephrase inequality constraints \mbox{$g^{[i]}({\xi}_{\cdot\mid k}^{[i]}) \leq 0$} as equality constraints $G({\xi}_{\cdot\mid k}^{[i]}) = 0$ (\cite{Sopasakis2020}) where for each constraint $\iota=1,\ldots,n_\text{ineq}$ holds
\begin{align}
G_\iota({\xi}_{\cdot\mid k}^{[i]}) \triangleq [g_\iota^{[i]}({\xi}_{\cdot\mid k}^{[i]})]_+.
\label{eq:DOCP_ADMM_PANOC_ineqPlusOp}
\end{align}
As \eqref{eq:DOCP_ADMM_PANOC_ineqPlusOp} is a nonsmooth function, we can not formulate $G({\xi}_{\cdot\mid k}^{[i]}) = 0$ as an ALM-type constraint (\cite{Sopasakis2020}). Instead, we apply the quadratic penalty method (PM) \cite[Chap.~17]{Nocedal2006} to replace the indicator function \eqref{eq:DOCP_decomp_indIneqCons} 
of the inequality constraints with the cost
\begin{align}
\phi_\text{PM}^{[i]}({\xi}_{\cdot\mid k}^{[i]};\alpha) {}\triangleq{} \frac{\alpha}{2} \, \| G({\xi}_{\cdot\mid k}^{[i]}) \|^2.
\label{eq:DOCP_ADMM_PANOC_PMcons}
\end{align}
With \eqref{eq:DOCP_ADMM_PANOC_ALMcons} and \eqref{eq:DOCP_ADMM_PANOC_PMcons}, we gain the Augmented Lagrangian function of the local OCP
\begin{align*}
\Phi^{[i]}({\xi}_{\cdot\mid k}^{[i]};\alpha,\mu) {}\triangleq{}     
\phi^{[i]}({\xi}_{\cdot\mid k}^{[i]}) + \phi_\text{ALM}^{[i]}({\xi}_{\cdot\mid k}^{[i]};\alpha,\mu) +\phi_\text{PM}^{[i]}({\xi}_{\cdot\mid k}^{[i]};\alpha) 
\end{align*}
which is a nonconvex $\mathcal{C}^1$ continuous differentiable function and as such a smooth reformulation of $\phi^{[i]}$. 
To compute a local solution of NLP \eqref{eq:DOCP_ADMM_algo_agentNLP}, we apply the open source code framework \texttt{OpEn v0.6.2} (\cite{Sopasakis2020}), available on \texttt{github.com/alphaville/optimization-engine}.
In an inner loop, \texttt{OpEn} utilizes the proximal averaged Newton method for optimal control (PANOC) to solve NLP \eqref{eq:DOCP_ADMM_algo_agentNLP} for a fixed penalty $\alpha$ and fixed Lagrangian multipliers $\mu$. In an outer loop, $\alpha$ and $\mu$ are updated to achieve constraint satisfaction as described in \cite{Sopasakis2020}. In every iteration of \prettyref{alg:DOCP_ADMM_algo} at time step $k$, the solver is warm-started with the solution of the previous iteration. Likewise, \texttt{OpEn} is warm-started at time step $k+1$ by exploiting the solution from time \mbox{step $k$.} For such kind of problems, PANOC has already shown superior performance (\cite{Stella2017,Katriniok2019a}).


\subsection{Convergence}
\label{sec:DOCP_convergence}

For convex consensus ADMM problems, that is, if every $\phi^{[i]}$ and $\gamma$ in \eqref{eq:DOCP_decomp_decoupNLP} are convex, convergence has been proven in literature (\cite{Boyd2011}). For nonconvex functions $\phi^{[i]}$ respectively $\Phi^{[i]}$, as in our case, we can show that the nonconvex consensus ADMM problem 
as well as the subproblems \eqref{eq:DOCP_ADMM_algo_agentNLP} 
converge to a set of stationary points under some mild conditions 
(\cite{Hong2016,Wang2019,Sopasakis2020}).
Without any assumptions on the iterates, \prettyref{alg:DOCP_ADMM_algo} is guaranteed to converge to a set of stationary points (i.e., to a local solution) if Problem \eqref{eq:DOCP_decomp_decoupNLP} meets certain regularity conditions and the step size $\rho$ 
is chosen large enough. 
By virtue of \cite{Hong2016}, our problem satisfies all assumptions which are a prerequisite to convergence, such as: the feasible set $({\Xi}_{k+j\mid k}^{[i]})_{i \in \mathcal{A},\,j \in \mathbb{N}_{[1,N]}}$ is closed and convex, $\gamma$ is convex and $\Phi^{[i]}$ has a Lipschitz-continuous gradient if 
the steering angle $\delta$ is constrained on the interval $(-\nicefrac{\pi}{2},\nicefrac{\pi}{2})$ --- technically, even tighter bounds are required, see \prettyref{sec:results_setup}. 
\section{Simulation Results}
\label{sec:results}

\subsection{Simulation Setup}
\label{sec:results_setup}

For a proof of concept, we evaluate the proposed consensus ADMM-based framework in a realistic lane change scenario with $N_A=3$ cooperative connected agents. According to \prettyref{fig:results_discussion_snapshots}a, the \SA{} (Agent 2, blue) is driving on the right (subject) lane while the \PAT{} (Agent 1, red) and \FAT{} (Agent 3, green) on the left (target) lane are initially blocking the road for a lane change of the \SA{}.  

In the simulation study, we adopt the control-oriented model \eqref{eq:modeling_stateSpace} as validation model. Every agent has dimensions $L^{[i]}=\unit[5]{m}$, $W^{[i]}=\unit[2]{m}$, \mbox{$l_f^{[i]}=\unit[1.4]{m}$} and $l_r^{[i]}=\unit[1.4]{m}$. The distance between the center of gravity and the front respectively rear bumper is $\unit[2.5]{m}$. In our scenario, we want the agents to keep a bumper-to-bumper distance of $\unit[10]{m}$. Taking vehicle dimensions and road curvature into account (the road radius is always larger or equal to $\unit[200]{m}$), we set $d_{\text{hw}}$ to $\unit[15]{m}$. This way, the minimum bumper-to-bumper distance is approximately $\unit[10]{m}$. For every agent, the initial and reference velocity is set to $\unitfrac[14]{m}{s}$. With a lane width of $w_\text{lane}=\unit[4]{m}$, the initial Frenet coordinates $(s^{[i]},\Delta y^{[i]})$ are $(\unit[12]{m},\unit[2]{m})$ and $(\unit[0]{m},\unit[2]{m})$ for the \PAT{} (red) and \FAT{} (green) in the left lane, respectively, while the initial coordinates of the \SA{} (blue) in the right lane are $(\unit[6]{m},\unit[-2]{m})$, see \prettyref{fig:results_discussion_snapshots}a.

In the consensus ADMM framework, we have selected the same weights for every agent: $q_{s}^{[i]} = 0$, $q_{\Delta y}^{[i]} = 1$, $q_{\Delta \psi}^{[i]} = 100$, $q_{v}^{[i]} = 1$, $r_{a_x}^{[i]} = 1$, $r_{\delta}^{[i]} = 600$ and $Q_N^{[i]}=Q^{[i]}$. The sample period between two consecutive runs of \prettyref{alg:DOCP_ADMM_algo} is set to $T_s=\unit[0.1]{s}$, the horizon length to $N=15$ and the penalty parameter to $\rho=100$. To keep the agents in their designated lanes, the bounds $( \Delta \underline{y}^{[i]}, \Delta \overline{y}^{[i]})$ on the lateral displacement are chosen as $(\unit[1.25]{m},\unit[2.75]{m})$ for the left lane and $(\unit[-1.25]{m},\unit[-2.75]{m})$ for the right lane. The absolute longitudinal acceleration should always be less or equal to $\unitfrac[4]{m}{s^2}$, the absolute lateral acceleration should not exceed 
$\unitfrac[3.5]{m}{s^2}$ and the total acceleration is bounded from above by $\unitfrac[4]{m}{s^2}$. Finally, the steering angle should be within $\pm \unit[5]{deg}$ while the maximum velocity is set to $\unitfrac[17]{m}{s}$. Simulations are run on an Intel i7 machine at $\unit[2.9]{GHz}$ with Matlab R2018b.


\subsection{Discussion of Results}
\label{sec:results_discussion}

\begin{figure}[t!] 
	\vspace*{-3mm}
	\hspace*{-3mm}
	\setlength\fwidth{0.42\textwidth}		
	\input{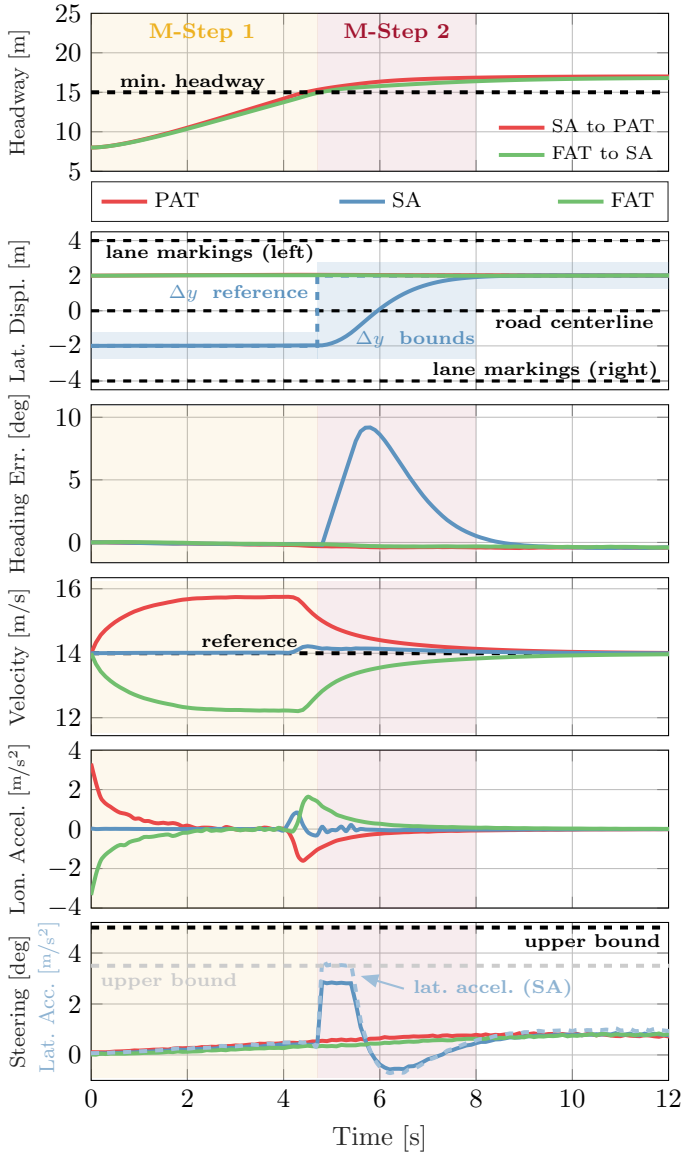}
	\vspace*{-8mm}
	\caption{Optimized state and input trajectories of the lane change maneuver. M-Step 1 is accomplished at \mbox{$t=\unit[4.7]{s}$} and M-Step 2 at $t=\unit[8.0]{s}$.} 
	\label{fig:results_discussion_statesInputs}
\end{figure}

\begin{figure*}[t] 
	\begin{center}
		\setlength\fwidth{0.80\textwidth}
		\input{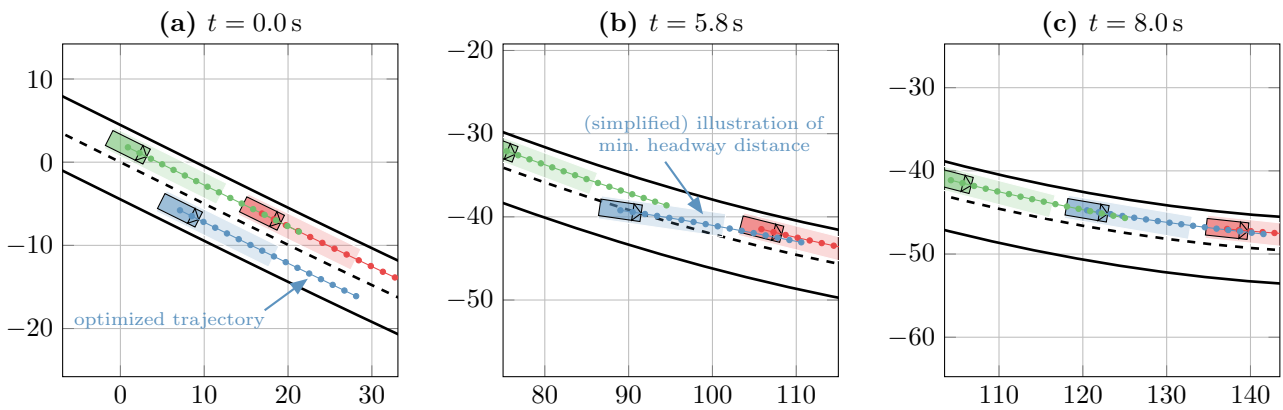}
		\vspace*{-2mm}
		\caption{Snapshots: 
			(Left) Initial configuration: \SA{} (blue) intends to change lanes;
			(Middle) \SA{} (blue) changes lanes after required headway has been established (M-Step 2);
			(Right) \SA{} (blue) has accomplished lane change maneuver.} 
		\label{fig:results_discussion_snapshots}
	\end{center}
\end{figure*}

In \prettyref{fig:results_discussion_statesInputs}, we illustrate the optimized state and input trajectories of the three agents while the steering angle plot (bottom plot, \prettyref{fig:results_discussion_statesInputs}) is augmented with the \SA{}'s lateral acceleration and its upper bound. Moreover, \prettyref{fig:results_discussion_snapshots} highlights three snapshots of the lane change maneuver.

\begin{figure}[t!] 
	\hspace*{-3mm}
	\setlength\fwidth{0.42\textwidth}
	\input{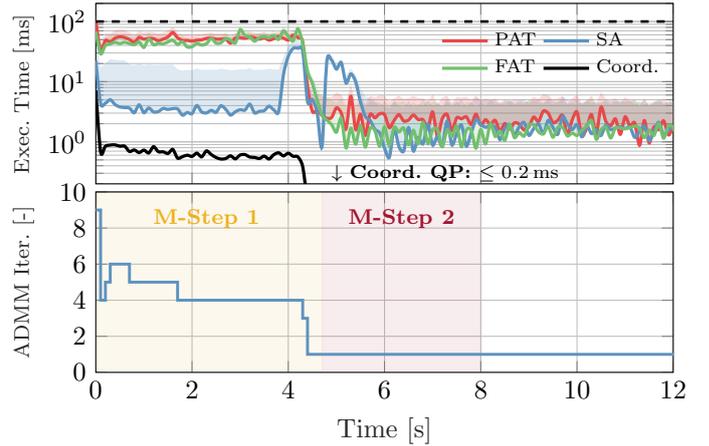}
	\vspace*{-8mm}
	\caption{Execution times (top) for solving the OCPs (solid lines) as a sum over all ADMM iterations (bottom) along with communication overhead (shaded area). The solving time of the coordination QP has been added to the \SA{} communication overhead (blue shaded area).} 
	\label{fig:results_discussion_compTime}
\end{figure}
During M-Step 1 of the maneuver, that is, until $t=\unit[4.7]{s}$ (light orange patch in \prettyref{fig:results_discussion_statesInputs}, \nth{1} plot) the agents establish the required headway distance to allow the \SA{} (blue) to change lanes safely. As a minimum cost maneuver, the \PAT{} (red) accelerates and the \FAT{} (green) decelerates to increase the headway distance whereas the \SA{} (blue) keeps its speed almost constant. As the \SA{} (blue) is located in the middle of the \PAT{} (red) and \FAT{} (green), the acceleration and speed trajectories of \PAT{} and \FAT{} are symmetric to each other --- which appears to be the most reasonable cooperative maneuver.
The maximum absolute longitudinal acceleration during M-Step 1 is $\unitfrac[3.3]{m}{s^2}$. 

At the beginning of M-Step 2 (light purple patch in \prettyref{fig:results_discussion_statesInputs}, \nth{1} plot), the \SA{} (blue) modifies its reference value $\Delta y^{\text{ref}}$ (dashed blue line in \prettyref{fig:results_discussion_statesInputs}, \nth{2} plot) from $\unit[-2]{m}$ (right lane) to $\unit[2]{m}$ (left lane). At the same time, the constraint bounds on the lateral displacement are adjusted to $(\Delta \underline{y}^{[i]},\Delta \overline{y}^{[i]})=(\unit[-2.75]{m},\unit[2.75]{m})$, see blue patch in \prettyref{fig:results_discussion_statesInputs}, \nth{2} plot. This enables the \SA{} (blue) to drive in both lanes and to initiate the lane change, see \prettyref{fig:results_discussion_snapshots}b. The lane change maneuver is accomplished after $\unit[3.3]{s}$, that is, at $t=\unit[8]{s}$, see also \prettyref{fig:results_discussion_snapshots}c. As we are rather focusing on lane change maneuvers than on platooning, there is no need to bound the headway distance from above to obtain a high vehicle density in the target lane. Consequently, the headway distance between the \SA{} (blue) and the \PAT{} (red) is always larger than the lower bound of $\unit[15]{m}$. 

Right after the end of M-Step 2, the agents' speed trajectory converge to the common set point of $\unitfrac[14]{m}{s}$.
As the maneuver is performed on a curved road, the wheel steering angle is almost always different from zero. During the entire maneuver, the agents' state and input trajectories are always smooth and within their designated bounds --- thus, satisfying all our requirements. 



Finally, the real-time capability of proposed control scheme, according to requirement R3 (see \prettyref{sec:DOCP}), needs to be assessed. \prettyref{fig:results_discussion_compTime} provides an overview of the agents' and the coordinator's execution times (top) and the required ADMM iterations (bottom). The execution times include the time needed to solve the local NLPs (solid lines in top plot) as well as the communication overhead related to the ADMM iterations (shaded area in same color, top plot). For the \SA{} (blue), the shaded area additionally incorporates the solving time of the coordination QP. As communication overhead, we assume a maximum end-to-end latency of \unit[3]{ms} as it has been specified by the 3rd Generation Partnership Project (3GPP) for the new 5G communication standard in order to support advanced driving functions (\cite{3gpp2019b}). 

\prettyref{fig:results_discussion_compTime} provides evidence that our control scheme is real-time capable as the maximum execution time (including communication overhead) is always below the sampling time of $\unit[100]{ms}$ (black dashed line). Only in the first time step, i.e., when the distributed problem is initialized, the computation of a first (initial) solution takes longer than the sampling time, that is, $\unit[131]{ms}$ for the \PAT{} (red). Initially, 9 ADMM iterations are required for convergence. However, the first time step can be viewed as an initialization phase which provides a viable initial guess for the second time step (in a receding horizon manner). During the rest of M-Step 1, we require 1 to 6 ADMM iterations for convergence while execution times for the \PAT{} (red) and \FAT{} (green) are mostly in a range of \unit[50]{ms} to \unit[70]{ms} --- with a maximum of \unit[87.5]{ms} at $t=\unit[4.2]{s}$ for the \FAT{} (green). The \SA{} (blue) shows lower execution times of \unit[15]{ms} to \unit[50]{ms}. During M-Step 2, consensus can be accomplished in a single ADMM step as the original minimum headway distance constraint \eqref{eq:COCP_centralizedNLP_CA} is always satisfied. Thus, computation times are much lower, that is, these are in a range of \unit[4]{ms} to \unit[7]{ms}. Only the \SA{} (blue) temporarily requires \unit[49]{ms} to solve its local OCP when the lateral acceleration constraints become active. 
\vspace*{-1mm}	
\section{Conclusion and Future Work}
\label{sec:conclusion}
\vspace*{-1mm}
We have introduced a distributed control concept for fully automated lane change maneuvers of CAVs which is able to accommodate situations when the traffic is dense and the target lane is already occupied. To this end, we have embedded the nonconvex problem formulation in a consensus ADMM framework which has shown convincing performance and real-time capability in our simulation study. As part of future work, we intend to remove the central coordinator from the control scheme, to validate the scheme in experiments and to embed it into a broader motion planning framework.



\bibliography{ifacconf}             

\begin{thebibliography}{19}
\providecommand{\natexlab}[1]{#1}
\providecommand{\url}[1]{\texttt{#1}}
\providecommand{\urlprefix}{URL }
\expandafter\ifx\csname urlstyle\endcsname\relax
  \providecommand{\doi}[1]{doi:\discretionary{}{}{}#1}\else
  \providecommand{\doi}{doi:\discretionary{}{}{}\begingroup
  \urlstyle{rm}\Url}\fi

\bibitem[{3GPP.org(2019)}]{3gpp2019b}
3GPP.org (2019).
\newblock Enhancement of {3GPP} support for {V2X} scenarios. {3GPP TS 22.186,
  V16.2.0}.
\newblock \urlprefix\url{https://portal.3gpp.org /desktopmodules/
  Specifications/SpecificationDetails.aspx ?specificationId=3180}.

\bibitem[{An and Jung(2018)}]{An2018a}
An, H. and Jung, J. (2018).
\newblock Design of a cooperative lane change protocol for a connected and
  automated vehicle based on an estimation of the communication delay.
\newblock \emph{Sensors}, 18(10).

\bibitem[{Bertsekas and Tsitsiklis(1989)}]{Bertsekas1989}
Bertsekas, D. and Tsitsiklis, J. (1989).
\newblock \emph{Parallel and Distributed Computation: Numerical Methods}.
\newblock {Prentice Hall}.

\bibitem[{{Bevly} et~al.(2016){Bevly}, {Cao}, {Gordon}, {Ozbilgin}, {Kari},
  {Nelson}, {Woodruff}, {Barth}, {Murray}, {Kurt}, {Redmill}, and
  {Ozguner}}]{Bevly2016a}
{Bevly}, D., {Cao}, X., {Gordon}, M., {Ozbilgin}, G., {Kari}, D., {Nelson}, B.,
  {Woodruff}, J., {Barth}, M., {Murray}, C., {Kurt}, A., {Redmill}, K., and
  {Ozguner}, U. (2016).
\newblock {Lane Change and Merge Maneuvers for Connected and Automated
  Vehicles: A Survey}.
\newblock \emph{IEEE Transactions on Intelligent Vehicles}, 1(1), 105--120.

\bibitem[{Blasi et~al.(2018)Blasi, K\"ogel, and Findeisen}]{Blasi2018a}
Blasi, S., K\"ogel, M., and Findeisen, R. (2018).
\newblock {Distributed Model Predictive Control Using Cooperative Contract
  Options}.
\newblock \emph{IFAC Conference on Nonlinear Model Predictive Control}, 51(20),
  448--454.

\bibitem[{Bock and Plitt(1984)}]{Bock1984}
Bock, H. and Plitt, K. (1984).
\newblock {A Multiple Shooting Algorithm for Direct Solution of Optimal Control
  Problems}.
\newblock \emph{IFAC World Congress}, 17(2), 1603--1608.

\bibitem[{Boyd et~al.(2011)Boyd, Parikh, Chu, Peleato, and Eckstein}]{Boyd2011}
Boyd, S., Parikh, N., Chu, E., Peleato, B., and Eckstein, J. (2011).
\newblock {Distributed Optimization and Statistical Learning via the
  Alternating Direction Method of Multipliers}.
\newblock \emph{Found. and Trends in Machine Learning}, 3(1), 1--122.

\bibitem[{Hong et~al.(2016)Hong, Luo, and Razaviyayn}]{Hong2016}
Hong, M., Luo, Z.Q., and Razaviyayn, M. (2016).
\newblock {Convergence Analysis of Alternating Direction Method of Multipliers
  for a Family of Nonconvex Problems}.
\newblock In \emph{SIAM Journal on Optimization}, volume~26, 337--364.

\bibitem[{Hu and Sun(2019)}]{Hu2019a}
Hu, X. and Sun, J. (2019).
\newblock Trajectory optimization of connected and autonomous vehicles at a
  multilane freeway merging area.
\newblock \emph{Transportation Research Part C: Emerging Technologies}, 101,
  111--125.

\bibitem[{Katriniok et~al.(2019)Katriniok, Sopasakis, Schuurmans, and
  Patrinos}]{Katriniok2019a}
Katriniok, A., Sopasakis, P., Schuurmans, M., and Patrinos, P. (2019).
\newblock {Nonlinear Model Predictive Control for Distributed Motion Planning
  in Road Intersections Using PANOC}.
\newblock In \emph{IEEE Conference on Decision and Control}, 5272--5278.

\bibitem[{Liu et~al.(2017)Liu, Ozguner, and Zhang}]{Liu2017a}
Liu, P., Ozguner, U., and Zhang, Y. (2017).
\newblock {Distributed MPC for cooperative highway driving and energy-economy
  validation via microscopic simulations}.
\newblock \emph{Transportation Research Part C: Emerging Techn.}, 77, 80--95.

\bibitem[{Nocedal and Wright(2006)}]{Nocedal2006}
Nocedal, J. and Wright, S.J. (2006).
\newblock \emph{Numerical Optimization}.
\newblock Springer, 2nd edition.

\bibitem[{{Qian} et~al.(2016){Qian}, {de La Fortelle}, and
  {Moutarde}}]{Qian2016}
{Qian}, X., {de La Fortelle}, A., and {Moutarde}, F. (2016).
\newblock {A hierarchical Model Predictive Control framework for on-road
  formation control of autonomous vehicles}.
\newblock In \emph{IEEE Intelligent Vehicles Symposium}, 376--381.

\bibitem[{Rajamani(2012)}]{Rajamani2012}
Rajamani, R. (2012).
\newblock \emph{{Vehicle Dynamics and Control}}, volume~2.
\newblock Springer.

\bibitem[{Sopasakis et~al.(2020)Sopasakis, Fresk, and Patrinos}]{Sopasakis2020}
Sopasakis, P., Fresk, E., and Patrinos, P. (2020).
\newblock {OpEn: Code Generation for Embedded Nonconvex Optimization}.
\newblock In \emph{IFAC World Congress}.

\bibitem[{Stella et~al.(2017)Stella, Themelis, Sopasakis, and
  Patrinos}]{Stella2017}
Stella, L., Themelis, A., Sopasakis, P., and Patrinos, P. (2017).
\newblock A simple and efficient algorithm for nonlinear model predictive
  control.
\newblock In \emph{IEEE Conference on Decision and Control}, 1939--1944.

\bibitem[{Wang et~al.(2016)Wang, Hu, Wang, Wang, Qin, and Bian}]{Wang2016}
Wang, D., Hu, M., Wang, Y., Wang, J., Qin, H., and Bian, Y. (2016).
\newblock Model predictive control-based cooperative lane change strategy for
  improving traffic flow.
\newblock \emph{Advances in Mechanical Engineering}, 8(2), 1--17.

\bibitem[{Wang et~al.(2019)Wang, Yin, and Zeng}]{Wang2019}
Wang, Y., Yin, W., and Zeng, J. (2019).
\newblock {Global Convergence of ADMM in Nonconvex Nonsmooth Optimization}.
\newblock \emph{Journal of Scientific Computing}, 78(1), 29--63.

\bibitem[{Wang et~al.(2017)Wang, Wu, and Barth}]{Wang2017a}
Wang, Z., Wu, G., and Barth, M. (2017).
\newblock {Developing a Distributed Consensus-Based Cooperative Adaptive Cruise
  Control System for Heterogeneous Vehicles with Predecessor Following
  Topology}.
\newblock \emph{Journal of Advanced Transportation}.

\end{thebibliography}

\end{document}